\documentclass[a4paper,12pt]{article}

\usepackage[T2A]{fontenc}
\usepackage[utf8]{inputenc}   % older versions of ucs package
\usepackage[russian]{babel}
\usepackage{marginnote}

\usepackage{CJKutf8}

\usepackage{amsmath,amstext,amsfonts,amsthm,amssymb}
\usepackage{eucal}
\usepackage{graphicx}
\usepackage{tikz-cd}
\usetikzlibrary{patterns}
\usetikzlibrary{calc}
\renewcommand{\subsubsection}[1]{\addtocounter{subsubsection}{1}
{\ \\[3pt]\bf \thesubsubsection. \  #1} }
\swapnumbers

{  \theoremstyle{definition}

}
\newcommand{\Spec}{\operatorname{Spec}}

\newcommand{\Tot}{\operatorname{Tot}}

\newcommand{\hra}{\hookrightarrow}

\newcommand{\lra}{\longrightarrow}

\newcommand{\dpar}{\partial}

\newcommand{\bea}{\begin{eqnarray*}}
\newcommand{\eea}{\end{eqnarray*}}
\newcommand{\bean}{\begin{eqnarray}}
\newcommand{\eean}{\end{eqnarray}}

\newcommand{\CA}{\mathcal{A}}

\newcommand{\CC}{\mathcal{C}}

\newcommand{\CD}{\mathcal{D}}
\newcommand{\CE}{\mathcal{E}}

\newcommand{\CL}{\mathcal{L}}
\newcommand{\CM}{\mathcal{M}}

\newcommand{\CO}{\mathcal{O}}
\newcommand{\CP}{\mathcal{P}}

\newcommand{\CS}{\mathcal{S}}

\newcommand{\CU}{\mathcal{U}}
\newcommand{\CV}{\mathcal{V}}

\newcommand{\CY}{\mathcal{Y}}
\newcommand{\CZ}{\mathcal{Z}}

\newcommand{\BC}{\mathbb{C}}

\newcommand{\BG}{\mathbb{G}}

\newcommand{\BZ}{\mathbb{Z}}

\newcommand{\nc}{\newcommand}

\nc{\Id}{\text{Id}}
\nc{\la}{\lambda}

\begin{document}

\centerline{\bf HOMOTOPY CHIRAL ALGEBRAS }

\bigskip\bigskip

\centerline{Fyodor Malikov and Vadim Schechtman }

\bigskip\bigskip

\begin{CJK}{UTF8}{min}

%\centerline{ヘッケの縮退代数とコ骸骨\footnote{hecke no shukutai daisu to kogaikotsu}}
%\centerline{ニルヘッケの代数とコ骸骨\footnote{niruhecke no daisu to kogaikotsu}}

\end{CJK}

\bigskip\bigskip

 \centerline{August 26, 2024}

\bigskip

\

%\hspace{8cm}{\it Гора родит смешную мышь }

\

%\hspace{7cm}{\it Parturiunt montes, nascetur ridiculus mus}

%\hspace{12cm}{\it Horatius}

%\hspace{4.5cm}{\it Пышутся горы родить, а смешной родится мышонок}

%\hspace{12cm}{\it Тредиаковский}

\

\hspace{1cm}{\it Here sat the man who agitated the scientific world with his Theory of Tittlebats} 

\

% это продолжение моего сегодняшнего письма; I will keep the notations from it. 

%{\it Чепуха совершенная делается на свете. Иногда вовсе нет никакого правдоподобия: вдруг тот самый нос, %который разъезжал в чине статского советника и наделал столько шуму в городе, очутился как ни в чем не %бывало вновь на своем месте, то есть именно между двух щек майора Ковалева.} 

\bigskip
\bigskip

%\begin{tikzpicture}[scale=.4, baseline=(current bounding box.center)]
 
% \node at (0,0){$\bullet$};
%  \draw (0,-4) -- (0,4); 
%  \draw (4,-0) -- (-4,0); 
 % \draw (4,4) --(4, -4) -- (-4,-4) -- (-4,4) -- (4,4); 
%  \node at (-.5, -0.5) {$0$}; 
% \node at (3.5, 3.5){$\hen$};  
%  \end{tikzpicture}
%  \quad\quad 
% { \huge  $\buildrel z^2\over  \lra$}
   \quad\quad

\centerline{\bf \S 1. Introduction} 

\

\

{\bf 1.1.} The aim of this paper is

(a)  to define a notion of a {\it homotopy chiral algebra} (HCA), or {\it $Ch_\infty$-algebra}, which means   
a chiral algebra up to higher homotopies;
 
(b)  to prove that the Cech complex $C^\bullet(\CU; \CA)$ of a sheaf $\CA$ 
of chiral algebras over a smooth algebraic variety $X$ admits a structure of a HCA. 

In particular the cohomology $H^\bullet(X; \CA)$ is a chiral algebra. 

\

{\bf 1.2.} Our idea is as follows. According to [BD] a chiral algebra  is a sort of a Lie algebra 
(cf. 1.3 below for the details). We know what a homotopy Lie algebra is (cf. [HS2] for example, and 1.3 below); 
adopting the definition to the chiral world, we get 1.1 (a). 

Now, the statement 1.1 (b) is the analog of the main result of [HS2] which says that the Cech complex 
of a sheaf of Lie algebras is a homotopy Lie algebra.  

{\bf 1.3.} Let us explain some details.
We will work over a ground field $k$ of characteristic $0$; $Mod(k)$ denotes the category of $k$-vector spaces. 

%We will deal with chiral algebras over a fixed 
%smooth algebraic curve $C$.

 The assertion 1.1 (b) is a consequence of a more general result.

 Let $\CL ie$ denote the Lie 
$k$-operad, so that algebras over $\CL ie$ are Lie algebras over $k$.

Let $\CC(k)$ denote the tensor category of complexes of $k$-modules. By definition 
if $A^\bullet, B^\bullet$ are two complexes, their tensor product $A^\bullet\otimes B^\bullet$ is a complex with components 
$$
(A^\bullet\otimes B^\bullet)^n = \oplus_{p+q = n}\ A^p\otimes B^q.
$$  

A {\it homotopy Lie operad} 
is an operad $\CL = \{\CL(n)\}$ in $\CC(k)$ together with a map of operads 
$\epsilon:\ \CL \lra \CL ie$ (where we consider each $\CL ie(n)$ as a complex concentrated in degree $0$) 
such that for each $n$
$$
\epsilon(n): \ \CL(n) \lra \CL ie(n)
$$
is a quasiisomorphism. By definition a {\it homotopy Lie algebra} is an algebra over a homotopy Lie operad.

The notion of a {\it pseudo-tensor category} generalizes that of an operad, cf. [BD] 1.1.1 and 2.1 below. 
A pseudo-tensor category with a unique object is the same as an operad.

Let $\CC$ be a pseudo-tensor $k$-category, cf. [BD] 1.1.7.   It is clear what a complex in $\CC$ is, and we denote by $\CC_{\CC(k)}$
the category of complexes in $\CC$, see 2.2 below. 

%\marginnote{ I've removed the equality 
%$\CC_{\CC(k)} = \CC(k)\otimes_k \CC.$; it is not clear to me how to make sense of it ... and I am not sure it is needed.}

Let $A^\bullet = \{ A^i\}$ be a cosimplicial Lie algebra in $\CC$. Let $M A^\bullet$ be the corresponding 
complex in $\CC$, i.e. 
$$
(MA)^i = A^i,\ i\geq 0
$$
and the differential is the alternating sum of faces. Thus $M A^\bullet$ is an object of $\CC_{\CC(k)}$. 

\

{\bf 1.4. Theorem.} {\it  $MA^\bullet$ admits a structure of a homotopy Lie algebra .}

This is Theorem 3.1.
It means that there exists a homotopy Lie operad $\epsilon:\ \CL \lra \CL ie$ such that 
$MA^\bullet$ admits a structure of an $\CL$-algebra. 

Applying this theorem to the Cech complex  $C^\bullet(\CU; \CA)$ of a sheaf of chiral algebras $\CA$ 
we get 1.1 (b). 

To prove 1.4  we follow the strategy of [HS2]. The main ingredient is a Lie version of the {\it Eilenberg - Zilber operad}, cf. 2.3 below.  See 3.1 for the details. 

The reader may find amusing the brief discussion of the {\it Eilenberg - Zilber PROP} in 2.7; 
it governs commutative and cocommutative bialgebras.
 
\

{\bf 1.5. Acknowledgement.} We are grateful to V.Hinich for helpful comments.

\

\newpage

\centerline{\bf \S 2. Recollections on pseudo-tensor categories etc.}

\

\

{\bf 2.1. Pseudo-tensor categories.} Recall (see [BD], 1.1) that a {\it pseudo-tensor $k$-category} $\CC$ consists of a class of objects 
$Ob \CC$ and the data below:

a rule that assigns to each finite non-empty set $I$, an $I$-tuple $\{ y_i\}_{i\in I}, y_i\in Ob \CC$, and $x\in Ob\CC$
a $k$-vector space
$$
P_I(\{y_i\}, x);
$$

a rule that assigns to each epimorphic map $\pi: J\lra I$ of finite sets and a family $\{ z_j\}_{j\in J}$ a $k$-polylinear composition map
$$
P_I(\{y_i\}, x)\times \Pi_I\ P_{J_i}(\{ z_j\}, y_i) \lra P_J(\{ z_j\}, x),
$$
$J_i := \pi^{-1}(i)$. 

These data should satisfy the axioms from [BD], 1.1.1. 

In particular (considering singletons $I = \{ * \}$) $\CC$ is a $k$-linear category, and the spaces 
$P_I(\{y_i\}, x)$ are functorial with respect to $y_i$ and $x$.

\

{\bf 2.2. The category $\CC_{\CC(k)}$}.  If $\CC$ is a pseudo-tensor $k$-category, then its objects are $k$-modules (use the spaces $P_1(x,x)$, $x\in \CC$); therefore one can consider complexes in $\CC$. Let $\CC_{\CC(k)}$ be the category of complexes
in $\CC$. 

For a finite set $I$ and a family of complexes $y^\bullet_i, x^\bullet$ the space of operations 
$$
P_I(\{y_i^\bullet\}, x^\bullet)
$$
is by definition the simple complex associated to the multiple complex
$$
P_I(\{y_i^{a_i}\}, x^b),\ a_i, b\in \BZ.
$$   

\

{\it Eilenberg - Zilber operad.}

\

{\bf 2.3.} Cf. [HS2], \S 5, [S]. 

$\Delta$ will denote the category of finite totally ordered sets, and $\Delta^o$ the opposite category. 

If $\CC$ is a category, $\Delta\CC$ will denote 
the category of cosimplicial objects in $\CC$, that is the category of functors $\Delta \lra \CC$. 

If $\CC$ is a tensor category, $\Delta\CC$ will be 
a tensor category, with $(X^\bullet\otimes Y^\bullet)^n = X^n\otimes Y^n$. 

Similarly $\Delta^o\CC$ will denote 
the category of simplicial objects in $\CC$, that is the category of functors $\Delta^o \lra \CC$.

Denote $\Delta(k) := \Delta Mod(k)$.  
Similarly $\Delta^o(k) := \Delta^o Mod(k)$; both these guys are tensor $k$-categories. 

We have a functor
$$
M^o:\ \Delta^o(k) \lra \CC^{\leq 0}(k),\ 
M^o(X)^p = X_{-p},\ 
$$
with 
$$
d: \ X_{-p}\lra X_{-p+1}
$$
being the alternating sum of faces.  (The change of sign ensures that the degree of the differential is +1, the convention we stick to throughout.)

For $n\geq 0$ $\Delta^n\in \Delta^o \CS ets$ will denote the standard simplex, and $k\Delta^n\in \Delta^o(k)$ 
its $k$-chains; denote
$$
C_\bullet(\Delta^n) := M^o(k\Delta^n)\in C^{\leq 0}(k).
$$
Varying $n$ we get an object  $Z \in \Delta \CC^{\leq 0}(k) = \CC^{\leq 0}(\Delta(k))$,
$$
Z([n]) = M^o(k\Delta^n).
$$
Here $[n] = \{0, 1, \ldots, n\}$.  

\

 On general grounds, $Hom_{\Delta(k)}(Z, A)$ is a $k$-vector space, but $Z$ being a simplicial complex it carries an endomorphism $d\in Hom_{\Delta(k)}(Z,Z)$
of square 0. Hence $Hom_{\Delta(k)}(Z, A)$ is a complex.

\

{\bf 2.3.1. Key lemma.} {\it For $A\in \Delta(k)$ we have
$$
M(A) = Hom_{\Delta(k)}(Z, A)\in C^{\geq 0}(k).
$$}

In other words, $Z$ represents the functor $M$. 

{\bf 2.3.2. Proof.}  Along with the functor $k\Delta^n\in \Delta^o(k)$, $k\Delta^n_m=k Hom_\Delta([m],[n])$, consider the represented functor $k\Delta_m\in\Delta(k)$,
$k\Delta^n_m=k Hom_\Delta([m],[n])$. We see that naturally $k\Delta_m\subset  Z$, and in fact $Z$ is a direct sum of  $k\Delta_m$, $m=1, 2, 3,...$. Therefore
$Hom_{\Delta(k)}(Z, A)\in C^{\geq 0}(k)$ is a direct product of  $Hom_{\Delta(k)}(k\Delta_m, A)$. By the Yoneda lemma, there is a vector space isomorphism
$$
A^m\stackrel{\sim}{\rightarrow}Hom_{\Delta(k)}(k\Delta_m, A);
$$
the lemma follows.  The proof is constructive, and an element $a_m\in A^m$ determines  the functor  morphism $Z\rightarrow A$ such that if 
$f\in Hom_{\Delta}([m],[n])\subset Z^n$, then $f$ is mapped to $A(f)a_m\in A^n$, and if $f\in Hom_{\Delta}([m'],[n])\subset Z^n$ with $m'\neq m$, then $f$ is mapped to 0.

$\square$  

\

{\bf 2.4. Eilenberg - Zilber operad.} More generally, for each $n\geq 0$ we may consider the tensor powers 
$Z^{\otimes n} \in \Delta \CC^{\leq 0}(k)$; to make sense of it we remark that $\CC^{\leq 0}(k)$ is a tensor category, and so $Z^{\otimes n}([m])=Z([m])^{\otimes n}$.

Define
$$
\CZ(n) = \Tot Hom_{\Delta(k)}(Z, Z^{\otimes n})\in \CC(k)
$$
where $\Tot$ means the total complex. 

These complexes form a $\CC(k)$-operad $\CZ$, to be called the {\it Eilenberg - Zilber operad}. 

Let $\CC om$ denote the commutative $k$-operad, with $\CC om(n) = k$ for all $n\geq 0$, so 
$\CC om$-agebras are commutative $k$-algebras.

%\newpage

Let $A$ be a cosimplicial commutative $k$-algebra. For each $n$ we have maps 
$$
Hom_{\Delta(k)}(Z, Z^{\otimes n})\otimes (MA)^{\otimes n} = 
Hom_{\Delta(k)}(Z, Z^{\otimes n})\otimes Hom_{\Delta(k)}(Z, A)^{\otimes n} \lra  
$$
$$
\lra Hom_{\Delta(k)}(Z, Z^{\otimes n})\otimes Hom_{\Delta(k)}(Z^{\otimes n}, A^{\otimes n}) \lra 
$$
$$
\lra Hom_{\Delta(k)}(Z, A^{\otimes n}) = M(A^{\otimes n}) 
\lra MA.
\eqno{(2.4.1)} 
$$
The last arrow is induced by the multiplication $A^{\otimes n} \lra A$. 
This defines a structure of a $\CZ$-algebra on $MA$. 

\ 

{\bf 2.5. } According to [HS1] or [HS2], 5.2.9, for each $n$ 
$$
H^i(\CZ(n)) = 0\ \text{if}\ i\neq 0,\ H^0(\CZ(n)) = k.
$$
In other words, $\CZ$ is a {\it homotopy $\CC om$ operad}  (to be precise, it is upon truncation, see 2.6.2 for details), and for each 
cosimplicial commutative $k$-algebra $A$ the complex $MA$ is a homotopy commutative algebra.

{\bf 2.5.1. Historical remark.} The operad $\CZ$ (without the name) was introduced and studied 
in the 1961 article by A.Dold, [D].

\

{\bf 2.6.} The Lie Eilenberg - Zilber operad 
$\CY = \{ \CY(n)\}$ is defined by 
$$
\CY(n) = \CZ(n)\otimes_k \CL ie(n). 
$$
By 2.5 $\CY$ is a homotopy Lie operad.

{\bf 2.6.1. Claim.} {\it Let $B$ be a cosimplicial Lie algebra. Then $MB$ admits a canonical $\CY$-action.}

{\bf Proof} cf. [HS2], 5.2.5, 5.2.6. We are given an action 
$$
\CL ie(n)\otimes B^{\otimes n} \lra B
$$
The required action
$$
\CY(n)\otimes MB^{\otimes n} \lra MB
$$
is defined as the following composition
$$
\CY(n)\otimes MB^{\otimes n} = Hom_\Delta(Z, Z^{\otimes n})\otimes \CL ie(n) \otimes 
Hom_\Delta(Z, B)^{\otimes n} \lra 
$$
$$
\lra Hom_\Delta(Z, Z^{\otimes n})\otimes 
Hom_\Delta(Z^{\otimes n}, \CL ie(n)\otimes B^{\otimes n}) \lra
$$
$$ 
\lra Hom_\Delta(Z, \CL ie(n)\otimes B^{\otimes n}) \lra Hom_\Delta(Z, B) = MB.
\eqno{(2.6.1.1)}
$$
It is a Lie analogue of (2.4.1).

\

{\bf 2.6.2. A truncation.} Let us define a suboperad
$$
\tau_{\leq 0}\CY = \{ \tau_{\leq 0}\CY(n)\}\subset \CY
$$
For each $n$ the inclusion
$$
\tau_{\leq 0}\CY(n)\hra \CY(n)
$$
is a quasiisomorphism.

On the other hand we have a canonical augmentation
$$
\epsilon(n):\ \tau_{\leq 0}\CY(n)\lra \CL ie(n)
$$
which a quasiisomorphism as well.  Thus $\tau_{\leq 0}\CY$ is a homotopy Lie operad.

\

{\bf 2.6.3. Corollary.} {\it  For each cosimplicial Lie algebra $B$ the complex 
$MB$ is an algebra for the operad $\tau_{\leq 0}\CY$, and therefore  
a homotopy Lie algebra.}

\

{\bf 2.7. Bialgebras.} More generally,  the complexes 
$$
\CZ(n,m) = \Tot Hom_{\Delta(k)}(Z^{\otimes m}, Z^{\otimes n})
$$
have only one, the $0$-th, nontrivial cohomology which is one-dimensional. 

They form a PROP 
$$
\CP\CZ = \{ \CZ(n,m),\ n, m\geq 0\}.
$$ 
in $C(k)$, cf. [A], 2.3, to be called the {\it Eilenberg - Zilber PROP}. 

Recall that a PROP encodes all possible operations in  bialgebras.  

Let $B$ be a cosimplicial commutative and cocommutative $k$-bialgebra. By the same token as in (2.4.1) 
we see that the complex $MB$ carries an action of $\CP\CZ$, so it may be called a {\it homotopy commutative and cocommutative bialgebra}.

\

\

%\newpage

\centerline{\bf \S 3. Finale: allegro con brio. The main theorem}

\
 
\

{\bf 3.1.} Let $\CC$ be a pseudo-tensor $k$-category. A $\CL ie$ algebra in $\CC$ is a pseudo-tensor category functor $F: \CL ie \rightarrow \CC$. In other words, it is an object $B\in \CC$ along with a collection of maps $F_n:\CL ie_n\rightarrow P_n(\{B\},B)$, $n=1,2,...$, which respect compositions.  This gives a category of $\CL ie$ algebras in $\CC$. 

 A cosimplicial Lie algebra
$B^\bullet$ is a functor from $\Delta$ to the  category of $\CL ie$ algebras in $\CC$, $[n]\mapsto B^n$. As before, we have the complex
 $MB^\bullet\in \CC_{\CC(k)}$ but there no identification $MB= Hom_{\Delta(k)}(Z, B)$ as  $Hom_{\Delta(k)}(Z, B)$  does not make sense, not immediately at least as $B$ is not a $k$-module (nor is it even a set.)  
 
 % \marginnote{I believe this is wrong; has been corrected in hca.3.1}
 
 Similarly, a homotopy  $\CL ie$ algebra in $ \CC_{\CC(k)}$ is a homotopy Lie operad $\CL$ and a pseudo-tensor category functor 
 $F:\CL\rightarrow  \CC_{\CC(k)}$. It amounts to an object $B\in  \CC_{\CC(k)}$ and a collection of maps $F_n:\CL(n)\rightarrow P_n(\{B\},B)$ which respect compositions.

{\bf  Theorem.} {\it Let $\CC$ be a pseudo-tensor $k$-category, $B^\bullet$ a cosimplicial $\CL ie$ algebra 
in $\CC$. Then the complex $MB^\bullet$ is a homotopy $\CL ie$ algebra in $\CC_{\CC(k)}$. More explicitly, there is a functor 
$\tau_{\leq 0}\CY \rightarrow  \CC_{\CC(k)}$, $\tau_{\leq 0}\CY(n)\rightarrow P_n(\{MB^\bullet\},MB^\bullet)$, $n=1, 2, 3, ...$.}
%\marginnote{Vadik, saying "it is an algebra over $\CY$" seems not quite correct as  this inadvertently hints at a tensor rather than a pseudo-tensor structure.}

{\bf 3.1.1. Beginning of the proof.} We will adjust the reasoning of 2.3 - 2.6 to the situation at hand where a tensor category structure is replaced with a  pseudo-tensor category structure.
More specifically, we will  rewrite  map (2.6.1.1)  
$$
\CY(n)\otimes MB^{\otimes n} \lra MB
$$
as a map
$$
\CY(n) \lra Hom_{\CC(k)}(MB^{\otimes n}, MB),
$$
thus obtaining an explicit formula for an operad morphism 
$$
\CY\lra \CE nd_{MB},
$$
where $\CE nd_{MB}$ is the endomorphism operad, $\CE nd_{MB}(n)=Hom_{\CC(k)}(MB^{\otimes n}, MB)$. We will see that this formula, albeit not the reasoning that led to it, makes sense for $\CE nd_{MB}(n)$ replaced with $P_n(\{MB\}, MB)$   and defines an operad morphism 
$$
\CY\lra P_\bullet(\{MB\}, MB).
$$
Here is this reasoning. We want to compute the composition map recorded in (2.6.1.1). By Key Lemma 2.3.1, or rather its proof, an element $b_m\in M^mB$ determines a morphism
$Z^k\rightarrow B^k$, for each $k$,  which sends $f:[m]\rightarrow [k]$ regarded as an element of $Z^k$ to $B(f)b_m\in B^k$. The tensor product   $b_{m_1}\otimes\cdot\otimes b_{m_n}\in (MB)^{\otimes n}$ is identified with the tensor product
of such morphisms, one for each $i$, $1\leq i\leq n$, which sends $f_i:[m_i]\rightarrow [k]$ to $B(f_i)b_{m_i}$. 

The first arrow in (2.6.1.1) makes this tensor product   into a morphism $Z^{\otimes n}\rightarrow B^{\otimes n}$, which sends, for each $k$,
\[
(Z^k)^{\otimes n}\ni f_{m_1}\otimes f_{m_2}\otimes\cdots\otimes f_{m_n}\mapsto B(f_{m_1})b_{m_1}\otimes\cdots\otimes  B(f_{m_n})b_{m_n}
\in (B^k)^{\otimes n}.
\eqno{(3.1.1.1)}
\]
In order to compute the next arrow in (2.6.1.1) we have to precompose morphism (3.1.1) with a morphism $\phi: Z\rightarrow Z^{\otimes n}$. 

By Key Lemma 2.3.1,
any such $\phi$ is a linear combination of those that are determined by its value $\phi(id_{[m]})$ for some fixed $m$. Assume that 
\[
\phi(id_{[m]})=\sum_{p_1,p_2,...,p_n}\alpha_{p_1,p_2,...,p_n}p_1\otimes\cdots\otimes p_n,
\eqno{(3.1.1.2)}
\]
for some $p_i: [l_i]\rightarrow m$ such that $l_1+l_2+\cdots+l_n=m$, $\alpha_{p_1,p_2,...,p_n}\in k$ are scalars.

The composition of (3.1.1.2) and (3.1.1.1) is the morphism $Z^m\rightarrow (Z^m)^{\otimes n}$ which sends

\[
id_{[m]}\mapsto \sum_{p_1,p_2,...,p_n}\alpha_{p_1,p_2,...,p_n}B(p_1)b_{m_1}\otimes\cdots\otimes B(p_n)b_{m_n}.
\eqno{(3.1.1.3)}
\]
The result might be nonzero only if $p_i\in Hom_\Delta([m_i],[m])$, as follows from 2.3.2, in which case, $B(p_i)b_{m_i}\in B^m$ for all $i$.

The last arrow in (2.6.1.1), essentially the action map $\CL ie(n)\otimes (B^m)^{\otimes n}\rightarrow B^m$ composed with (3.1.1.3), defines a morphism
$Z\rightarrow B$ determined by the assignment
\[
id_{[m]}\mapsto \sum_{p_1,p_2,...,p_n}\alpha_{p_1,p_2,...,p_n}\psi(B(p_1)b_{m_1}\otimes\cdots\otimes B(p_n)b_{m_n});
\eqno{(3.1.1.4)}
\]
here $\psi\in \CL ie(n)$ and $\psi(q_1,...,q_n)\in B^m$ is the value of the action map $\CL ie(n)\otimes (B^m)^{\otimes n}\rightarrow B^m$  on
$\psi\otimes(q_1\otimes\cdots\otimes q_n)$.

Here is the operadic meaning of (3.1.1.4). 
Let us give ourselves an arbitrary $\psi\in\CL ie(n)$ and  $\phi\in Hom_{\Delta(k)}(Z,Z^{\otimes n})$. We can assume that $\phi$, such as (3.1.1.2), defines  an element 
\[
\bar{\phi}\in\otimes_{i=1}^n Hom_k(B^{m_i}, B^m)
\]
 by virtue of the cosimplicial structure of  $B^\bullet$. (Simply because $B(p_i)\in Hom_k(B^{m_i},B^{m})$.)
Because of this, formula (3.1.1.4) says that the operad morphism 
\[
\CL ie\otimes\CY\rightarrow \CE nd_{MB},
\]
 whose existence is asserted by Claim 2.6.1, 
sends $\psi\otimes\phi\in \CL ie(n)\otimes Hom_{\Delta(k)}(Z,Z^{\otimes n})$ to an element of 
\[
\CE nd_{MB}(n)=\prod_{l_1,...,l_n;l_{n+1}} Hom_k(\otimes_{i=1}^{n}B^{l_i},B^{l_{n+1}})
\]
such that the entry in $Hom_k(\otimes_{i=1}^{n}B^{l_i},B^{l_{n+1}})$ is 0 unless, $l_i=m_i$, $l_{n+1}=m_1+\cdots+ m_n$, $i=1,...,n$, in which case it equals
\[
\psi_m(\bar{\phi});
\eqno{(3.1.1.5)}
\]
here $\psi_m(\bar{\phi})$ is the value of the composition
\[
Hom_k(\otimes_{i=1}^n B^m, B^m)\otimes\otimes_{i=1}^n Hom_k(B^{m_i},B^m)\rightarrow Hom_k(\otimes_{i=1}^n B^{m_i}, B^m)
\eqno{(3.1.1.6)}
\]
on $\psi_m\otimes\bar{\phi}$, where $\psi_m\in Hom_k(\otimes_{i=1}^n B^m, B^m)$ is the image of $\psi\in\CL ie(n)$ under the action (operad) morphism
$\CL ie\rightarrow \CE nd_{B^m}$ arising by virtue of the Lie algebra structure on $B^m$.

It remains to notice that formulas (3.1.1.5-6), unlike their origin, (2.6.1.1), or  their derivation,   are written purely in terms of a pseudo-tensor structure on the category 
$\CC(k)$, where
$P_I(\{L_i\},M)=Hom_{\CC(k)}(\otimes_I L_i,M)$. This allows us to upgrade them to the case of a cosimplicial Lie object in an arbitrary pseudo-tensor category.

{\bf 3.1.2. Refinement of the statement.}  Let us return to the setting of Theorem 3.1, where $B^\bullet$ is a cosimplicial Lie algebra in an arbitrary pseudo-tensor category $\CC$. We define a family of maps, one for each $n=1,2,3,...$,
\[
F_n:\CL ie(n)\otimes\CY(n) \rightarrow  P_n(\{MB\},MB)=\prod_{l_1,...,l_n;l_{n+1}}P_n(\{B^{l_1},...,B^{l_n}\}, B^{l_{n+1}})
\]
so that
for each $\psi \in \CL ie(n)$ and $\phi\in Hom_{\Delta(k)}(Z,Z^{\otimes n})$ the entry of $F_n(\psi\otimes \phi)$ in $P_n(\{B^{l_1},...,B^{l_n}\}, B^{l_{n+1}})$ equals
\[
F_n(\psi\otimes \phi)_{l_1,...,l_n;l_{n+1}}=\left\{
\begin{array}{ccc}
\psi_m(\bar{\phi})&\mbox{ if }& l_i=m_i, l_{n+1}=m_1+\cdots m_n, 1\leq i\leq n\\
0&\mbox{ otherwise; }&
\end{array}
\right.
\eqno{(3.1.2.1)}
\]
here $\psi_m$ is the image of $\psi\in \CL ie(n)$ under the structure map $\CL ie(n)\rightarrow P_n(\{B^m\}, B^m)$, 
$
\bar{\phi}\in\otimes_{i=1}^n Hom(B^{m_i}, B^m)
$ is the element determined by $\phi$
as in 3.1.1, especially (3.1.1.5), and $\psi_m(\bar{\phi})$ is the value of the
pseudo-tensor structure composition map (cf. the beginning of 2.1)
\[
P_n(\{B^m\}, B^m)\otimes \otimes_{i=1}^n Hom_\BC(B^{m_i},B^m)\rightarrow P_n(\{B^{m_1}, B^{m_2},...,B^{m_n}\}, B^m)
\]
on the element $\psi_m\otimes\bar{\phi}$.

The refinement of Theorem 3.1 is as follows:

{\bf 3.1.3. Theorem.} {\it The collection of maps $\{F_n,\; n=1,2,...\}$ defines a pseudo-tensor functor $\CL ie\otimes \CY\rightarrow \CC_{\CC(k)}$, hence a homotopy Lie
algebra structure on $MB$.} 

{\bf Proof }  follows from (3.1.2.1) in straightforward manner. 

\

{\bf 3.2.} Let $C$ be a smooth curve over $k$, $\CC = \CM(C)$ the category of $\CD$-modules over $C$. 
It carries two pseudo-tensor structures: one defined in [BD], 22.3, and another one defined 
in [BD], 3.1. $Lie$ objects of the former one are called $Lie^*$-algebras, those of the latter - chiral algebras over $C$.

 Let $B$ is a cosimplicial $Lie^*$ (resp. chiral) algebra over $C$. By 3.1 the complex $MB$ is a homotopy 
$Lie^*$ (resp. chiral) algebra over $C$.

\

{\bf 3.3.} Let $X = \bigcup\ U_i$, $\CU = \{ U_i \}$, be an open covering of a smooth algebraic variety 
$X$ over $k$, and $\CA$ a sheaf of chiral $C$-algebras over $X$. By definition the Cech complex 
$$
C^\bullet(\CU; \CA) = M\tilde{C}^\bullet(\CU; \CA)
$$
where $\tilde{C}^\bullet(\CU; \CA)$ 
is a cosimplicial chiral $C$-algebra with
$$
\tilde{C}^p(\CU; \CA) = \oplus_{i_1,\ldots,i_p} \Gamma(U_{i_1}\cap\ldots\cap U_{i_p}; \CA).
$$   

It follows from 3.2 that  $C^\bullet(\CU; \CA)$ is a homotopy chiral $C$-algebra.

\

\

%\newpage

\centerline{\bf \S 4. Epilogue. Higher Borcherds operations}

\

\

{\bf 4.1. Jacobi identity up to a homotopy.} Let $V^\bullet$ be a homotopy Lie algebra in $\CC(k)$.

By definition this means that  
we are given an operad $\CL = \{\CL(n)\}$ in the tensor category $\CC(k)$  
together with a map of operads 
$$
\epsilon: \CL \lra \CL ie, 
$$
such that for all $n$ the map $\epsilon(n): \CL(n) \lra \CL(n)$ is a quasi-isomorphism (here we consider 
$\CL ie(n)$ as complexes concentrated in degree $0$), 
and the complex $V^\bullet$ is an $\CL$-algebra, that is,  we are given a map of operads
%\marginnote{!}
\[
F:\CL\rightarrow \CE nd_{V^\bullet} \eqno{(4.1.0)}
\]
\[
F_n: \CL(n)\rightarrow Hom_{\CC(k)}((V^\bullet)^{\otimes n},V)
\]

Recall that we can identify $\CL ie(n)$ with the space of Lie polynomials in variables $e_1, \ldots, e_n$ 
containing exactly one letter $e_i$ for each $i$, cf. [BD], 1.1.10.  

Thus we have a unique up to  coboundary  $0$-cocycle $c\in \CL(2)^0$
such that $\epsilon(c) = [e_1, e_2]$.  
Pick such a cocycle. It gives rise to a bracket (a skew symmetric map), cf. (4.1.0),
%\marginnote{!}
$$
[ , ]\stackrel{\mbox{def}}{=}F_2(c):\ V^\bullet\otimes V^\bullet \lra V^\bullet.
$$ 
Let $j = j_c\in \CL(3)^0$ be the element corresponding to the Lie polynomial  
$$
[e_1,[e_2,e_3]] + [e_2,[e_3,e_1]] + [e_3,[e_1,e_2]].
\eqno{(4.1.1)} 
$$
All three elements corresponding to the summands of (4.1.1) are cocycles since operadic compositions 
in $\CL$ are maps of complexes. 

Then $\epsilon(j) = 0$ and $dj = 0$. 
Therefore there exists $j'\in \CL(3)^{-1}$ such that $dj' = j$.    

Define
%\marginnote{!} 
$J=F_3(j)\in Hom(V^{\bullet\otimes 3}, V)$, i.e. an arrow 
$$
J:\ V^{\bullet\otimes 3}\lra V^\bullet,
$$
which corresponds to $j$ via (4.1.0), i.e.
$$
J(a, b, c) = [a, [b, c]] + [b, [c, a]] + [c, [a, b]]
$$
(a "jacobiator"). 
 
Let  similarly
$$
J':\ V^{\bullet\otimes 3}\lra V^\bullet [-1]
$$
be the map corresponding to $j'$. 

Then
%\marginnote{!} 
by definition
$$
J = d\circ J' + J'\circ d,
$$
in other words 
$$
d J'(a, b, c) + J'(da, b, c) + (-1)^{\deg a} J'(a, db, c) + (-1)^{\deg a + \deg b}  J'(a, b, dc) = 
$$
$$
= [a, [b, c]] + [b, [c, a]] + [c, [a, b]].
$$

\

(A) {\it $Lie^*$ algebras case}  

\

{\bf 4.2.} Now let us write down the analogue of 4.1 with the tensor category $\CC(k)$ replaced by 
the category of $\CM(C)$ of right $\CD$-modules over a curve $C$ equipped with the $Lie^*$ or the chiral pseudo-tensor structure. 

We will start with a slightly simpler case of $Lie^*$-algebras (see [BD], 0.3). 

Recall that for a finite 
nonempty set $I$ and a family $M_i, i\in I, N$ in $\CM(C)$ the space of $*$-operations 
$$
P_I^*(\{ M_i\}, N) = Hom_{\CM(C^I)}(\boxtimes M_i, \Delta_*^{(I)}N)
$$
where $\Delta^{(I)}:\ C\hra C^I$ is the embedding of the principal diagonal. 

Let $C$ be the affine line $\Spec k[t]$, and let $T = \BG_a$ denote the group of translations which operates on $C$.

Let 
%\marginnote{!} 
$\CV$ be a translationally invariant $Lie^*$ algebra over 
$C$, so $\CV$  is a right $\CD_C$-module. 

Let $V^\ell = V\omega_C^{-1}$ denote the corresponding left $\CD_C$-module. Let  $V$ be 
the space of $T$-invariant  global sections $V = \Gamma(C;\CV^\ell)^T$, or otherwise the stalk $\CV^\ell_x$  at any  $x\in C$, cf. [BD], 0.15.

Consider a commutative $k$-subalgebra 
$$
R = k[\dpar]\subset \Gamma(C,\CD_C) = k[t,\dpar].
$$
It is a cocommutative Hopf algebra with a comultiplication
$$
R\lra R\otimes_k R,\ \dpar \mapsto \dpar\otimes 1 + 1\otimes\dpar.
$$
The space $V$ is an $R$-module, and 
for a finite nonempty set $I$ we have
$$
P_I^*(\{ \CV\}, \CV) = Hom_{R^{\otimes I}}(V^{\otimes I}, V\otimes_R R^{\otimes I}),
\eqno{(4.2.1)}
$$
the composition of $*$ operations being induced by the Hopf algebra structure on $R$, cf. [BD], 1.1.13 and 
1.1.14.

We denote
$$
R^{\otimes I} = k[\dpar_1,\ldots, \dpar_n],
$$
and we identify 
$$
V\otimes_R R^{\otimes I} = V\otimes_k k[\dpar_1,\ldots,\dpar_{n-1}]
$$
by sending
$$
x\otimes \dpar_n \mapsto x\otimes_R(\dpar - \sum_{i=1}^{n-1} \dpar_i) = x\dpar\otimes 1 - 
\sum_{i=1}^{n-1} x\otimes\dpar_i.
\eqno{(4.2.2)}
$$
 
For example the $*$-bracket on $\CV$ is the same as an arrow 
$$
\{ , \}:\ V\otimes V \lra V\otimes_R R^{\otimes 2}. 
$$
For $a, b\in V$ we write, by making use of the agreement (4.2.2), 
$$
\{a , b\} = \sum_{n\geq 0} a_{(n)}b\otimes \dpar_1^n/n!
\eqno{(4.2.3)}
$$  
which defines a set of  operations
$$
_{(n)}:\ V\otimes V \lra V,\ n\in \BZ_{\geq 0}.
$$

\

{\bf 4.2.1. Caim.} {\it For any $a, b$ $a_{(n)}b = 0$ for sufficiently big $n$.}

Indeed, 
$$
V\otimes_R R^{\otimes 2} = V\otimes_k R,
$$
and $R$ is a polynomial ring, so the sum in the RHS of (4.2.3) must be finite. $\square$ 

\

{\bf 4.2.2.} The $*$-operations can be composed; as an example, we consider the following 3 elements of $Hom_{R^{\otimes 3}}(V^{\otimes 3},V\otimes_RR^{\otimes 3})$:
$$
a\otimes b\otimes c\mapsto \{a,\{b,c\}\}=\sum_{m,n\geq 0}a_{(m)}(b_{(n)}c)\otimes \dpar_1^m/m! \dpar_2^n/n!,
$$
$$
a\otimes b\otimes c\mapsto \{b,\{a,c\}\}=\sum_{m,n\geq 0}b_{(n)}(a_{(m)}c)\otimes \dpar_1^m/m! \dpar_2^n/n!,
$$
$$
a\otimes b\otimes c\mapsto \{\{a,b\},c\}\}=\sum_{p,q\geq 0}(a_{(p)}b){(q)}c)\otimes \dpar_1^p/p! (\dpar_1+\dpar_2)^{p+q}/(p+q)!.
$$
Collecting the like terms, the left hand side of the Jacobi identity can be written as the following element of $Hom_{R^{\otimes 3}}(V^{\otimes 3},V\otimes_RR^{\otimes 3})$:
\begin{eqnarray*}
a\otimes b\otimes c&\mapsto& \{a, \{b, c\}\} - \{b, \{a, c\}\} -\{\{a, b\}, c\}\\
&=& 
\sum_{m,n\geq 0}\left([a_{(m)}, b_{(n)}]c- \sum_{j\geq 0}\binom{m}{j} \left(a_{(j)}b\right)_{(m + n - j)}c\right)\otimes \dpar_1^m/m! \dpar_2^n/n!
\end{eqnarray*}

Equating with 0 the coefficient of $\dpar_1^m/m! \dpar_2^n/n!$ we see that the Jacobi identity 
$$
\{a, \{b, c\}\} - \{b, \{a, c\}\} = \{\{a, b\}, c\}
$$
is equivalent to a family  of  identities labelled by pairs of integers $(m,n)$:
$$
[a_{(m)}, b_{(n)}] c= \sum_{j\geq 0}\binom{m}{j} (a_{(j)} b)_{(m + n - j)}c,
\eqno{(4.2.4)} 
$$
cf. [K], (2.7.1).

Similarly, the skew symmetry 
$$
\{ b, a\} = - \{ a, b\}
$$
is translated to 
$$
a_{(n)}b = (-1)^{n+1}\sum_{j\geq 0} (b_{(n+j)}a)\dpar^j/j!.
\eqno{(4.2.5)}
$$
The sum here is finite due to Claim 4.2.1.

\

Thus $V$ is what was dubbed a {\it vertex Lie algebra} in [P], and a  {\it conformal algebra} in [K], Def. 2.7.

\

%{\bf 4.2.3. Claim.} {\it The sheaf $\CD iff_C$ of differential operators is a $Lie^*$ coalgebra.} 

%EXPLAIN WHAT DOES IT MEAN
  
%Indeed, we can identify $R$ with the space 
%$$
%\Gamma(C, \CD iff^T_C)
%$$
%where $\CD iff^T_C$ is the sheaf of translationally invariant differential operators. 

%$R$ carries a structure of a conformal algebra in the sense of [K], which is equivalent to saying that 
%$\CD iff_C$ is a $Lie^*$ algebra.

%CHECK IF THIS IS TRUE!

% Note that by construction, (4.2.3), $a_{(n)}b=0$ for all $a,b\in V$ if $n\gg 0$, which ensures that all %potentially infinite summations in the formulas written above are in fact finite.

{\bf 4.3.} Now let us pass to the homotopy picture. 
 
 Let $\CV^\bullet$ be a cosimplicial translationally invariant $Lie^*$  algebra on $C$; 
it gives rise to a cosimplicial object $V^\bullet = \Gamma(C; \CV^{\bullet\ell})^T$.   

Let 
$$
MV^\bullet = \Gamma(C; M\CV^{\bullet\ell})^T
$$
be  
the corresponding complex. By 3.1 $M\CV^\bullet$ will be a homotopy Lie algebra in $\CM(C)$, i.e. the complex $M\CV^\bullet$ will be equipped with a structure of 
an algebra  over some homotopy Lie operad $\CL$. 

As in 4.1 we pick a cocycle $c\in \CL(2)^0$ mapping 
to the Lie bracket in $\CL ie(2)$ under the structure map $\CL(2)\rightarrow \CL ie(2)$. It gives rise to an operation 
$$
[ , ] = [ , ]_c\in P_2(\{M\CV^\bullet, M\CV^\bullet\}, M\CV^\bullet)^0.
\eqno{(4.3.1)}
$$
The bracket (4.3.1) gives rise to a set of operations
$$
_{(n)}:\ MV^{\bullet\otimes 2} \lra MV^\bullet,
$$
$n\in \BZ_{\geq 0}$. 

These operations satisfy the identities (4.2.4)  
up to a homotopy. To explain them let us start with a Lie polynomial
$$
[e_1, [e_2, e_3]] - [e_2, [e_1, e_3]] - [[e_1, e_2], e_3] 
$$
Then acting as in 4.1 we get an  element $j'\in \CL(3)^{-1}$ which gives rise 
to an operation
$$
J'\in P_3(\{M\CV^\bullet, M\CV^\bullet, M\CV^\bullet\}, M\CV^\bullet)^{-1}).
$$
such that
$$
d J'(a, b, c) + J'(da, b, c) + (-1)^{\deg a} J'(a, db, c) + (-1)^{\deg a + \deg b}  J'(a, b, dc) = 
$$
$$
= \{a, \{b, c\}\}-   \{b, \{a, c\}\} + \{\{a,b\},c\}.
$$
Passing to fibers at some point $x\in C$ and using the conventions of (4.2.1-2) we get a set of  
{\it secondary operations}
$$
j'_{(m n)}:\ (MV^\bullet)^{\otimes 3} \lra MV^\bullet[-1],\ 
m, n\in \BZ_{\geq 0}
$$ 
such that
$$
J'(a,b,c)=\sum_{m,n\geq 0}j'(a,b,c)_{(mn)}\otimes \dpar_1^m/m! \dpar_2^n/n!.
$$
Operations $j'_{(mn)}$ must satisfy
$$
d j'(a, b, c)_{(m n)} + j'(da, b, c)_{(m n)} + (-1)^{\deg a} j'(a, db, c)_{(m n)} + (-1)^{\deg a + \deg b}  j'(a, b, dc)_{(m n)} = 
$$
$$
= a_{(m)}(b_{(n)}c) - b_{(n)}(a_{(m)}c) - \sum_{j\geq 0}\binom{m}{j}(a_{(j)}b)_{(m+n-j)}c
\eqno{(4.3.2)}
$$

\

(B) {\it Chiral algebras case}

\

{\bf 4.4. From chiral algebras to vertex algebras.} For a finite 
nonempty set $I$ and a family $M_i, i\in I, N$ in $\CM(C)$ the space of chiral operations 
$$
P_I^{ch}(\{ M_i\}, N) = Hom_{\CM(C^I)}(j_*j^*\boxtimes_I M_i, \Delta_*^{(I)}N)
$$
where $\Delta^{(I)}:\ C\hra C^I$ is the embedding of the principal diagonal and 
$$
j: C^I \setminus \bigcup \Delta_{ij}\hra C^I 
$$
is the open embedding of the complement to all diagonals.

Notice that the definition of the $D$-module direct image functor, 
$$ 
\Delta_*^{(I)}N=N\otimes_{D_C}\Delta^{(I)\ast}D_{C^I},
$$
implies
\[
 \Delta_*^{(I)}N=Diff(\CO_{C^I},N),
 \]
 the space of differential operators on $C^I$ with values in $N$, which is regarded as an $\CO_{C^I}$-module by pull-back. We have 
 \begin{eqnarray*}
P_I^{ch}(\{ M_i\}, N) &=& Hom_{\CM(C^I)}(j_*j^*\boxtimes_I M_i, \Delta_*^{(I)}N)\\
&=&
 Hom_{D_{C^I}}(\boxtimes_I M_i\otimes_{\CO_{C^I}}j_*j^*\CO_{C^I},Diff(\CO_{C^I},N))\\
 &\hookrightarrow& Hom_{D_{C^I}}(\boxtimes_I M_i, Hom_k(\CO_{C^I}\otimes_{\CO_{C^I}}j_*j^*\CO_{C^I},N))\\
 &=& Hom_{D_{C^I}}(\boxtimes_I M_i, Hom_k(j_*j^*\CO_{C^I},N))
\end{eqnarray*}

Let $C = \Spec k[t]$ be the affine line, and suppose that $M_i, N$ are translationally invariant; 
let $V_i, W$ be the corresponding spaces of translationally invariant sections. In this case, the result of the computation above is simplified to the effect that
$$
Hom_{\CM(C^I)}(j_*j^*\boxtimes_I M_i, \Delta_*^{(I)}N)\subset Hom_k((\otimes_I V_i), Hom_k(k[(t_i - t_j)^{-1}], W),
$$
or, equivalently,
$$
Hom_{\CM(C^I)}(j_*j^*\boxtimes_I M_i, \Delta_*^{(I)}N)\subset Hom_k((\otimes_I V_i)[(t_i - t_j)^{-1}], W).
$$

% EXPLAIN WHY -- {\bf ADJUNCTION!}

\

{\bf 4.4.1.} 
A chiral algebra is a Lie algebra object $M$ for this pseudo-tensor structure. Let $M$ be translationally invariant, and $V$ be the corresponding space of translationally invariant sections.  
 
The chiral bracket
$$
[ , ] \in P_2^{ch}(\{ M, M\}, M)  
$$
defines a map 
$$
\phi_{[,]}:\ (V\otimes V)[(t_1 - t_2)^{ \pm 1}] \lra V.
\eqno{(4.4.1)}
$$
For each $p\in \BZ$, $a, b\in V$ 
let us denote  
$$
a_{(p)}b := \phi_{[,]}((a\otimes b)(t_1 - t_2)^p)
$$ 

Let us introduce a generating function 
%\marginnote{!}
$$
a(t)b := \sum_{n = -\infty}^\infty \frac{a_{(n)}b}{t^{n+1}}.
\eqno{(4.4.2)}
$$
In other words the generating series is an operation: 
$$
\sum_{n\in\BZ} a_{(n)}b t^{-n-1}: V\lra V\{ t\} := \prod_{n = -\infty}^{\infty} Vt^n.
$$
Then $a_{(p)}b$ will be expressed as a residue 
$$
a_{(p)}b = \phi_{[,]}((a\otimes b)(t_1 - t_2)^p)= \int_{C(t_2)} (t_1 - t_2)^p a(t_1 - t_2) b dt_1,
$$
the contour integral being taken along a circle with center $t_2$.

{\bf 4.4.2. Jacobi identity and Borcherds relations.}
The iterated residues allow us  to compute 
compositions of brackets.

For example, the double bracket
$$
[ , [ , ] ] \in P_3^{ch}(\{ M, M, M\}, M)
$$
defines and is defined by a map
$$
\phi_{[ , [ , ] ]}:\ V^{\otimes 3}[(t_1 - t_2)^{\pm 1}, (t_1 - t_3)^{\pm 1}, (t_2 - t_3)^{\pm 1}] \lra V
$$
which is computed as a double residue.

 Namely, for each $(p, q, r)\in \BZ^3$
 %\marginnote{!}
$$
\phi_{[ , [ , ] ]}(a\otimes b\otimes c)((t_1 - t_2)^{p}(t_1 - t_3)^{q}(t_2 - t_3)^{r}) = 
$$ 
$$
\int_{C_{t_1}(t_3); |t_1 - t_3| > |t_2 - t_3|}\int_{C_{t_2}(t_3)} (t_1 - t_2)^p(t_1 - t_3)^{q}(t_2 - t_3)^{r}a(t_1 - t_3)b(t_2 - t_3)c dt_2 dt_1 =   
$$
$$
= \sum_{j\geq 0} \binom{p}{j}((-1)^j a_{(p+q-j)}b_{(r+j)}c
$$
%\marginnote{!} 
This is because in the indicated domain, $|t_1-t_3|>|t_2-t_3|$, we keep $t_1-t_3$, $t_2-t_3$, expand
\[
(t_1-t_2)^p=(t_1-t_3-(t_2-t_3))^p=\sum_{j=0}^\infty{p\choose j}(t_2-t_3)^j(t_1-t_3)^{p-j},
\]
and then use (4.4.2).

% CHECK! --{\bf CHECKED!!}

%\marginnote{!}  
Similarly with $a\leftrightarrow b$ and $t_1\leftrightarrow t_2$ simultaneously interchanged.

The third chiral double bracket
$$
\{\{ a, b\}, c\}\in P_3^{ch}(\{ M, M, M\}, M)
$$
gives rise to 
$$
\phi_{ [[ , ], ]}:\ V^{\otimes 3}[(t_1 - t_2)^{\pm 1}, (t_1 - t_3)^{\pm 1}, (t_2 - t_3)^{\pm 1}] \lra V 
$$
which is computed by the double residue:
%\marginnote{!}
$$
\phi_{ [[ , ], ]}(a\otimes b\otimes c)((t_1 - t_2)^{p}(t_1 - t_3)^{q}(t_2 - t_3)^{r}) = 
$$
$$
= \int\int_{|t_1 - t_2| < |t_2 - t_3|} (t_1 - t_2)^p(t_1 - t_3)^q(t_2 - t_3)^r(a(t_1 - t_2)b)(t_2 - t_3) c =  
$$
$$
= \sum_{j\geq 0} \binom{q}{j}(a_{(p+q-j)}b)_{(r+j)}c 
$$

%\marginnote{!} 
This is because this time around we keep $t_1-t_2$, $t_2-t_3$, expand
\[
(t_1-t_3)^q=(t_1-t_2+(t_2-t_3))^q=\sum_{j=0}^\infty{q\choose j}(t_2-t_3)^j(t_1-t_2)^{q-j},
\]
and then use (4.4.2).

Let us now consider the left-hand side of the Jacobi identity, that is, the following ternary chiral operation
$$
a\otimes b\otimes c\mapsto Jac(a,b,c)=[a, [b, c]] - [b, [a, c]] - [[a, b], c].
$$
Combining the above-written formulas we obtain\marginnote{!}
$$
\phi_{Jac}(a\otimes b\otimes c)((t_1 - t_2)^{p}(t_1 - t_3)^{q}(t_2 - t_3)^{r}) = 
$$
$$
\sum_{j\geq 0} \binom{p}{j}((-1)^j a_{(p+q-j)}b_{(r+j)}c - (-1)^{j+p} b_{(r+p-j)}a_{(q+j)}c) -
$$
$$
- \sum_{j\geq 0} \binom{q}{j}(a_{(p+q-j)}b)_{(r+j)}c
$$

We see that the Jacobi identity for the chiral bracket 
$$
[a, [b, c]] - [b, [a, c]] = [[a, b], c]
$$
is equivalent to the following family of identities parametrized by triples of integers $(p,q,r)$
$$
\phi_{Jac}(a\otimes b\otimes c)((t_1 - t_2)^{p}(t_1 - t_3)^{q}(t_2 - t_3)^{r}) = 0;
$$
therefore, as we have just computed, to the family of
 Borcherds identities 
 %\marginnote{!}
$$
\sum_{j\geq 0} \binom{p}{j}((-1)^j a_{(p+q-j)}b_{(r+j)}c - (-1)^{j+p} b_{(r+p-j)}a_{(q+j)}c) =
$$
$$
= \sum_{j\geq 0} \binom{q}{j}(a_{(p+q-j)}b)_{(r+j)}c
$$

%OPERADIC LANGUAGE: 

%$R$ AS A MEROMORPHIC  HOPF ALGEBRA? 

%COMPARE 4.2.3 

\

{\bf 4.5. Homotopy context: the secondary Borcherds operations.} Let $\CV^\bullet$ be a cosimplicial translationally invariant chiral algebra 
over $C = \Spec k[t]$. By Thm 3.1 the corresponding complex $M\CV^\bullet$ is a homotopy chiral algebra. 

Let $V^\bullet$ be the corresponding cosimplicial space of $T$-invariant global sections, and $MV^\bullet$ 
the corresponding complex.

Proceeding similarly to 4.3 we get a ternary chiral homotopy operation
$$
Jac'\in P_3^{ch}(\{ M\CV^\bullet, M\CV^\bullet, M\CV^\bullet\}, M\CV^\bullet [-1]), 
$$
a "secondary Jacoiator"\ , 
which gives rise to 
$$
\phi_{Jac'}:\ MV^{\bullet\otimes 3}[(t_1 - t_2)^{\pm 1}, (t_1 - t_3)^{\pm 1}, (t_2 - t_3)^{\pm 1}] 
\lra MV^\bullet[-1]. 
$$
For $a, b, c\in V^\bullet,\ m, k, l\in \BZ$ we denote 
%\marginnote{! - and below}
$$
Jac'(a\otimes b\otimes c)_{(p, q, r)} := \phi_{Jac'}((t_1 - t_2)^{p}(t_1 - t_3)^{q}(t_2 - t_3)^{r}),
$$  
so we get 
a collection of operations 
$$
_{(p,q,r)}:\ MV^{\bullet\otimes 3} \lra MV^\bullet[-1], 
$$
of degree $-1$, 
the "secondary Borcherds operations"\ .

For $a \in V^\alpha, b\in V^\beta, c\in V^\gamma$ we have    
$$
(d\circ _{(p,q,r)} + _{(p,q,r)}\circ d)(a\otimes b\otimes c) = 
$$
$$
\sum_{j\geq 0} \binom{p}{j}((-1)^j a_{(p+q-j)}b_{(r+j)}c - (-1)^{j+p+\alpha\beta} b_{(r+p-j)}a_{(q+j)}c) -\sum_{j\geq 0} \binom{q}{j}(a_{(p+q-j)}b)_{(r+j)}c
\eqno{(4.5.1)}
$$

These are the {\it secondary Borcherds identities}. 

\

The following assertion has been formulated in [B], without a detailed proof.

\

{\bf 4.6. Corollary.} {\it Let $X$ be a topological space (or a Grothendieck topology), and 
$\CV$ a sheaf of vertex algebras over $X$. Then the Cech cohomology 
$^cH^\bullet(X, \CV)$ is a $\BZ_{\geq 0}$-graded vertex algebra.}

\

\

\centerline{\bf References}

\bigskip\bigskip

[A] J.F.Adams, Infinite loop spaces, Princeton University Press and University of Tokyo Press, 
Princeton, NJ, 1978. 

[BD] A.Beilinson, V.Drinfeld, Chiral algebras, AMS Coll. Publ. {\bf 51}, AMS Providence, Rhode Island, 2004.

[B] L.Borisov, Vertex algebras and mirror symmetry, {\it Comm. Math. Phys., } {\bf 215} (2001), 517  - 557.

[D] A.Dold, \"Uber die Steenrodschen Kohomologieoperationen, {\it Ann. Math.} {\bf 73} (1961), 258 - 294. 

[HS1] V.Hinich, V.Schechtman, On homotopy limit of homotopy algebras, Yu.I.Manin (Ed.), $K$-theory, Arithmetic and Geometry, Moscow, USSR 1984 - 1986,  LNM {\bf 1289}, 240 - 264.

[HS2] V.Hinich, V.Schechtman, Homotopy Lie algebras, I.M.Gelfand Seminar, {\it Adv. Sov. Math.} {\bf 16}, 
Part 2 (1993), 1 - 28.

[K] V.Kac, Vertex algebras for beginners, University Lecture Series {\bf 10}, AMS Providence, Rhode Island, 
1997.

[M] F.Malikov, Strongly homotopy chiral algebroids, {\it Contemp. Math.} {\bf 711} (2018), AMS Providence, 
Rhode Island, 1 - 35. 

[P] M.Primc, Vertex algebras generated by Lie algebras, {\it J. Pure Appl. Alg.} {\bf 135}  (1999), 
253 - 293. 

[S] V.A.Smirnov, On the cochain complex of topological spaces, {\it Math. USSR Sbornik} 
{\bf 43} (1982) no. 1, 133  - 144. 

\

F.M.: Dept. of Mathematics, University of Southern California, Los Angeles CA 90089, USA

Email: {\tt fmalikov@usc.edu}

\smallskip

V.S.: Institut de Math\'ematiques de Toulouse, Universit\'e Paul Sabatier, 118 Route de Narbonne, 
31062 Toulouse, France; Kavli IPMU, 5-1-5 Kashiwanoha, Kashiwa, Chiba, 277-8583 Japan

Email: {\tt schechtman@math.ups-tlse.fr}

%\newpage

%\input ???

%  WRITE HERE

\end{document}